\documentclass[12pt]{article}
\usepackage{amsfonts, amsmath}
\author{Gleb G. Gusev\thanks{Partially supported by the grants
RFBR-007-00593, INTAS-05-7805, NWO-RFBR 047.011.2004.026,
RFBR-08-01-00110-a.}}
\title{Monodromy zeta-functions of deformations and Newton diagrams}
\date{}

\DeclareMathOperator{\pt}{pt}

\newtheorem{thm}{Theorem}

\newtheorem{prop}{Proposition}

\newenvironment{dfntn} {\smallskip\noindent{\bf
Definition\/}.}{\smallskip\par}

\newenvironment{exmpls} {\smallskip\noindent{\bf
Examples\/}.}{\smallskip\par}

\newenvironment{rmrks} {\smallskip\noindent{\bf Remarks\/}.}{\smallskip\par}
\newenvironment{prf} {\noindent{\em Proof\/}.}{{ $\Box$}\smallskip\par}
\newenvironment{acknl} {\smallskip\noindent{\bf Acknowledgement\/}.}{\smallskip\par}

\renewcommand{\a}{\alpha }

\renewcommand{\d}{\delta }
\newcommand{\D}{\Delta }
\newcommand{\e}{\varepsilon }

\newcommand{\G}{\Gamma }
\renewcommand{\l}{\lambda }
\renewcommand{\L}{\Lambda }

\newcommand{\s}{\sigma }
\newcommand{\Sig}{\Sigma}

\newcommand{\z}{\zeta }
\newcommand{\p}{\pi}

\newcommand{\xx}{\mathbf x}
\newcommand{\zz}{\mathbf z}
\newcommand{\yy}{\mathbf y}
\newcommand{\uu}{\mathbf u}

\newcommand{\kk}{\mathbf k}
\newcommand{\R}{\mathbb{R}}
\newcommand{\C}{\mathbb{C}}
\newcommand{\DD}{\mathbb{D}}

\newcommand{\Id}{\mathop{\mathrm{Id}}\nolimits}

\begin{document}
\maketitle

\abstract{For a one-parameter deformation of an analytic complex
function germ of several variables, there is defined its
monodromy zeta-function. We give a Varchenko type formula for this
zeta-function if the deformation is non-degenerate with respect to
its Newton diagram.}

\section{Introduction}

Let $F$ be the germ of an analytic function on $(\C^{n+1}, 0)$,
where $\,\C^{n+1} = \C_\s\times \C^n_{\zz}$, $\s$ is the
coordinate on $\C$, and $\,\zz = (z_1,z_2,\ldots, z_n)\,$ are the
coordinates on $\C^n$. The germ $F$ provides a deformation
$\,f_{\s} = F(\s, \cdot)\,$ of the function germ $\,f = f_0\,$ on
$(\C^n, 0)$. We give formulae for the monodromy zeta-functions of
the deformations of the hypersurface germs $\,\{f=0\}\cap
(\C^*)^n\,$ and $\,\{f = 0\}\,$ at the origin in terms of the
Newton diagram of $F$. A reason to study deformations of
hypersurface germs and their monodromy zeta-functions was inspired
by their connection with zeta-functions of deformations of
polynomials: \cite{GZ-2}.

Let $A$ be the complement to an arbitrary analytic hypersurface
$Y$ in $\C^n$: $\,A = \C^n\setminus Y$. Let $\,V\, = \,\{F=0\}\cap
(\C_\s\times A) \cap B_\e$, where $\,B_\e \subset \C^{n+1}\,$ is
the closed ball of radius $\e$ with the centre at the origin. Let
$\,\DD^*_\d \subset \C_\s\,$ be the punctured disk of radius $\d$
with the centre at the origin. For $\,0< \d \ll \e\,$ small enough
the restriction to $V$ of the projection $\,\C^{n+1} \to \C_\s\,$
onto the first factor provides a fibration over $\,\DD^*_\d\,$
(\cite{LeDT}). Denote by $V_c$ the fibre over the point $c$.
Consider the monodromy transformation $\,h_{F,A} \colon V_c \to
V_c\,$ of the above fibration restricted to the loop $\,c \cdot
\exp(2\p it),\, t\in [0,1],\, |c|$ is small enough.

The zeta-function of an arbitrary transformation $\,h\colon X\to
X\,$ of a topological space $X$ is the rational function
$\,\zeta_h(t) = \prod_{i\geq 0} (\det(\Id - th_*|_{H^c_i(X;
\C)}))^{(-1)^i}$, where $\,H^c_i(X; \C)\,$ is the $i$-th homology
group with closed support.

\begin{dfntn}
The zeta-function of the monodromy transformation $h_{F,A}$ will
be called the monodromy zeta-function of the deformation $f_\s$ on $A$:
$\zeta_{{f_\s}|_{A}}(t) = \zeta_{h_{F,A}}(t)$.
\end{dfntn}

For a power series $\,S= \sum c_\kk \yy^\kk$, $\,\yy^\kk =
y_1^{k_1}\cdots y_m^{k_m},\,$ one defines its Newton diagram as
follows. Denote by $\,\R_+ \subset \R\,$ the set of non-negative
real numbers. Denote by $\G'(S)$ the convex hull of the union
$\cup_{c_\kk \neq 0} (\kk + \R_+^m).$ The Newton diagram of the
series $S$ is the union of compact faces of $\G'(S)$. For a germ
$G$ on $\,\C^m\,$ at the origin, its Newton diagram $\G(G)$ is the
Newton diagram of its Taylor series at the origin.

For a generic germ $F$ on $\,(\C^{n+1}, 0)\,$ with fixed Newton
diagram $\G\in \R_+^{n+1}$ the zeta-functions
$\,\zeta_{{f_\s}|_{(\C^*)^n}}(t),\,\zeta_{{f_\s}|_{\C^n}}(t)\,$
are also fixed. We provide explicit formulas for these
zeta-functions in terms of the Newton diagram $\G$.

\section{The main result (a Varchenko type formula)}

Let $F$ be a germ of a function on $(\C^{n+1},0)$. Let $\,\kk =
(k_0, k_1, \ldots, k_n)\,$ be the coordinates on $\R^{n+1}$
corresponding to the variables $\s, z_1,\ldots, z_n$ respectively.
For $\,I\subset \{0, 1,\ldots, n\}\,$ denote by
$\,\R^I,\,\G^I(F)\,$ the sets $\,\{\kk \mid k_i=0,\, i\notin
I\}\subset \R^{n+1}\,$ and $\,\G(F)\cap \R^I\,$ respectively.

An integer covector is called primitive if it is not a multiple of
another integer covector. Let $P^I${\label{Z}} be the set of
primitive integer covectors in the dual space ${(\R^I)}^*$ such
that all their components are strictly positive. For $\,\a \in
P^I,\,$ let $\G^I_\a(F)$ be the subset of the diagram $\G^I(F)$
where $\a|_{\G^I(F)}$ reaches its minimal value: $\G^I_\a (F) =
\{\xx\in \G^I(F) \mid \a(\xx)=\min(\a|_{\G^I(F)})\}$ (for
$\,\G^I(F)=\emptyset\,$ we assume $\,\G^I_\a (F)=\emptyset$).
Consider the Taylor series of the germ $F$ at the origin: $\,F=
\sum F_{\kk} \s^{k_0} z_1^{k_1}\ldots z_n^{k_n}.\,$ Denote:
$\,F_\a = \sum_{\kk \in \G_\a^{\{0,1,\ldots ,n\}}} F_{\kk}
\s^{k_0} z_1^{k_1}\ldots z_n^{k_n}$.

\begin{dfntn}
A germ $F$ of a function on $\,(\C^{n+1},0)\,$ is called
non-degenerate with respect to its Newton diagram if for any $\,\a
\in P^I\,$ the $1$-form $\,dF_\a\,$ does not vanish on the germ
$\,\{F_\a = 0\}\cap (\C^*)^{n+1}\,$ at the origin (see
\cite{Var}).
\end{dfntn}

For $\,I\in \{0,1,\ldots, n\}\,$ such that $0\in I$, we denote:
$$
\zeta_F^I (t)= \prod_{\a\in P^I} (1-t^{\a(\frac{\partial}{\partial
k_0})})^{(-1)^{l-1}l!\,V_l(\G^I_\a (F))},
$$
where $\,l = |I|-1,\,$
$\frac{\partial}{\partial k_0}$ is the vector in $\R^I$ with the
single non-zero coordinate $\,k_0 = 1,\,$ and $V_l(\cdot)$ denotes
the $l$-dimensional integer volume, i.e. the volume in a rational
$l$-dimensional affine hyperplane of $\R^I$ normalized in such way
that the volume of the minimal parallelepiped with integer
vertices is equal to 1. We assume that $\,V_0(\pt)=1\,$ and for
$\,n\geq 0$ one has $\,V_n(\emptyset)=0$.

\begin{thm}\label{thm1}
Let $F$ be non-degenerate with respect to its Newton diagram
$\G(F)$. Then one has
\begin{equation}
\label{1}\zeta_{{f_\s}|_{(\C^*)^n}}(t) = \zeta_F^{\{0,1\ldots,
n\}}(t),
\end{equation}
\begin{equation}
\label{1a}\zeta_{{f_\s}|_{\C^n}}(t) = (1-t) \times \prod_{I\colon
0\in I\subset \{0,1, \ldots, n\}} \zeta_F^I (t).
\end{equation}
\end{thm}

\begin{rmrks}
1. The equation (\ref{1}) implies the equation (\ref{1a}) because
of the following multiplicative property of the zeta-function. Let
$\,h\colon X \to X\,$ be a transformation of a CW-complex $X$. Let
$\,Y\subset X\,$ be a subcomplex of $X$. Assume that
$\,h(Y)\subset Y,\,$ $h(X\setminus Y)\subset (X\setminus Y).\,$
Then $\zeta_{h|_X}(t)= \zeta_{h|_{X\setminus Y}}(t)\times
\zeta_{h|_{Y}}(t)$.

One can see that $\,\zeta_{f_\s|_{\{0\}}} (t) = (1-t)\times
\zeta_F^{\{0\}} (t)$. In fact, in the case $\,\G^{\{0\}} =
\emptyset\,$ one has $\,\zeta_{f_\s|_{\{0\}}}(t)= (1-t)\,$,
$\,\zeta_F^{\{0\}} (t) = 1.\,$ Otherwise $\,\zeta_{f_\s|_{\{0\}}}
(t) = 1,\,$ $\,\zeta_F^{\{0\}} (t) = (1-t)^{-1}$.

2. The zeta-function $\z_{f_\s|_{\C^n}}(t)$ coincides with the
monodromy zeta-function of the germ of the function $\s\colon
\{F=0\} \to \C_\s\,$ at the origin. The main theorem of \cite{Oka}
provides a formula for the zeta-functions of germs of functions on
complete intersections in non-degenerate cases. One can apply this
formula to the germ $\s$ and verify that the formula (\ref{1a})
agrees with the one of M. Oka. But (\ref{1a}) can not be deduced
from the result of M. Oka because the function $\s$ doesn't
satisfy the condition of "convenience" (\cite[p. 17]{Oka}).
\end{rmrks}

\begin{exmpls}
1. Let $\,F(\s,\zz) = f(\zz)-\s$. The monodromy zeta-function of
the deformation $f_\s$ coincides with the (ordinary) monodromy
zeta-function $\zeta_f(t)$ of the germ $f$ on $(\C^n,0)$ (see,
e.g., \cite{Var}). In this case the $\,l$-dimensional faces
$\,\G^I_\a (F)\,$ (where $\,l= |I|-1 > 0$) are cones of integer
height $1$ over the corresponding $\,(l-1)$-dimensional faces
$\G^{I\setminus \{0\}}_{\a|_{\{k_0=0\}}} (f)$. One has:
$$
V_l(\G^I_\a (F)) = V_{l-1}(\G^{I\setminus
\{0\}}_{\a|_{\{k_0=0\}}}(f))\,/\,l
$$
with $\,\a(\partial/\partial k_0) = \min(\a|_{\G^{I\setminus \{0\}}(f)})$. This means that in this case the equation (\ref{1a}) coincides with the Varchenko
formula (\cite{Var}).

2. For a deformation $\,F(\s, \zz)\,$ of the form $\,f_0(\zz)- \s
f_1(\zz)\,$ the fibre
$$
(\{\s\}\times  \{f_\s = 0\})\cap B_\e
$$
is the disjoint union of the sets
$$
(\{\s\}\times \{\,f_0/f_1 =
\s\})\cap B_\e
$$
and
$$
(\{\s\}\times \{f_0=f_1=0\})\cap B_\e.
$$
If $\,f_0(0) =  f_1(0) = 0,\,$ then $\,\zeta_{f_\s|_{\C^n}}(t)=
(1-t)\times \zeta_{(f_0/f_1)|_{\C^n}}(t),\,$ otherwise
$\,\zeta_{f_\s|_{\C^n}}(t)= \zeta_{(f_0/f_1)|_{\C^n}}(t)\,$ (the
zeta-function of the meromorphic function $f_0/f_1$: \cite{GZ}).

For $I\subset \{0,1,\ldots, n\}$ such that $\,0\in I,\,$ and for a
covector $\a\in P^I,\,$ assume that the face $\,\G^I_\a (F)\,$ has
dimension $l$, where $\,l= |I|-1>1$. Then $\,\G^I_\a (F)\,$ is the
convex hull of the corresponding faces $\,\D^I_{\a,0} =
\{0\}\times \G^{I\setminus \{0\}}_{\a|_{\{k_0=0\}}}(f_0)\,$ and
$\,\D^I_{\a,1} = \{1\}\times \G^{I\setminus
\{0\}}_{\a|_{\{k_0=0\}}}(f_1)$, which lie in the hyperplanes
$\,\{k_0=0\}\,$, $\,\{k_0=1\}\,$ respectively. It is not difficult
to show (see, e.g., \cite[Lemma 1]{G}) that $\,V_l(\G^I_\a (F)) =
V^I_\a /l$, where
\begin{multline*}
V^I_\a = V_{l-1}(\D^I_{\a,0}, \ldots,\D^I_{\a,0}) + V_{l-1}(\D^I_{\a,0}, \ldots, \D^I_{\a,0},
\D^I_{\a,1})+\\
+ \ldots + V_{l-1}(\D^I_{\a,0}, \D^I_{\a,1}, \ldots,
\D^I_{\a,1})+ V_{l-1}(\D^I_{\a,1}, \ldots,\D^I_{\a,1}).
\end{multline*}
Here $V_{l-1}$ denotes the $(l-1)$-dimensional Minkowski's mixed
volume: see, e.g., \cite{Oka}). Moreover, $\,\a(\partial/\partial
k_0) = \min(\a|_{\G^{I\setminus \{0\}}(f_0)})-
\min(\a|_{\G^{I\setminus \{0\}}(f_1)})$, thus (\ref{1a}) coincides
with the main result of \cite{GZ}.
\end{exmpls}

\section{A'Campo type formula}

Proof of Theorem \ref{thm1} uses an A'Campo type formula
(\cite{Campo}) written in terms of the integration with respect to
the Euler characteristic (\cite{GZ-2}).

For a constructible function $\Phi$ on a constructible set $Z$
with values in a (multiplicative) Abelian group $G$ its integral
$\,\int_Z \Phi^{d\chi}\,$ with respect to the Euler characteristic
$\chi$ is defined as $\prod_{g\in G} g^{\chi(\Phi^{-1}(g))}$ (see
\cite{Viro}). Further we consider $\,G = {\C(t)}^*\,$ to be the
multiplicative group of non-zero rational functions in the
variable $t$.

Let $F$ be a germ of an analytic function on $\,(\C^{n+1},0)$
defined on a neighbourhood $U$ of the origin. Let $Y$ be a
hypersurface in $\C^n$. Denote $S = (\C_\s\times Y) \cup
\{\s=0\}$. Consider a resolution $\,\p\colon (X,D)\to (U,0)\,$ of
the germ of the hypersurface $\,\{F=0\}\cup S\,$ at the origin,
where $D = \p^{-1}(0)$ is the exceptional divisor.

\begin{thm}\label{thm2}
Assume $\p$ to be an isomorphism outside of $\,\p^{-1}(U\cap S)$.
Then
\begin{equation}
\label{3} \zeta_{f_\s|_{\C^n\setminus Y}}(t) = \int_{D\cap W}
{\zeta_{\Sig|_{W\setminus Z},\,x}(t)}^{d\chi},
\end{equation}
where $W$ is the proper preimage of $\,\{F=0\}\,$ (i.e. the
closure of $\,\p^{-1}(V)$, $\,V = \left ((\{F=0\}\cap U) \setminus
S\right )$), $\Sig = \s \circ \p$, $Z = \p^{-1}(\C_\s\times Y)\,$
and $\,\zeta_{\Sig|_{W\setminus Z},\,x}(t)\,$ is the monodromy
zeta-function of the germ of the function $\Sig$ on the set
$\,W\setminus Z\,$ at the point $\,x\in D\cap W$.
\end{thm}

\begin{prf}
The map $\p$ provides an isomorphism $\,W\setminus (Z\cup \{\Sig =
0\}) \to V$, which is also an isomorphism of fibrations provided
by the maps $\Sig$ and $\s$ over sufficiently small punctured
neighbourhood of zero $\DD^*_\d \subset \C_\s$. Therefore the
monodromy zeta-functions of this fibrations coincide,
$\,\zeta_{f_\s|_{\C^n\setminus Y}}(t) = \zeta_{\Sig|_{W\setminus
Z}}(t)$ (the monodromy zeta-function of the "global" function
$\Sig$ on $W\setminus Z$).

Applying the localization principle (\cite{GZ-2}) to $\Sig$ we
obtain:
\begin{equation}
\label{4}\zeta_{f_\s|_{\C^n\setminus Y}}(t) = \int_{W\cap \{\Sig
=0\}}{\zeta_{\Sig|_{W\setminus Z},\,x}(t)}^{d\chi}.
\end{equation}

The integration is multiplicative with respect to subdivision of
its domain. One has $\,W\cap \{\Sig=0\} = (D\cap W) \sqcup((W\cap
\{\Sig =0\})\setminus D)$. Thus the right hand side of (\ref{4})
is the product $\left [\int_{D\cap W}{\zeta_{\Sig|_{W\setminus
Z},\,x}(t)}^{d\chi}\right ] \cdot \left [\int_{W\cap (\{\Sig
=0\}\setminus D)}{\zeta_{\Sig|_{W\setminus
Z},\,x}(t)}^{d\chi}\right ]$. The first factor coincide with the
right hand side of (\ref{3}); we prove that the second factor
equals 1.

For a point $\,x\in D,\,$ its neighbourhood $\,U(x)\subset X\,$
with a coordinate system $\,u_1,u_2,\ldots,u_{n+1}$ is called {\it
convenient} if each of manifolds $\,D, Z\,$ can be defined on
$U(x)$ by an equation of type $\,\uu^{\kk} = 0\,$ and each of
functions $\,\Sig, \,\tilde{F} = F\circ \p\,$ has the form
$\,a\,\uu^{\kk}$, where $a(0)\neq 0$. One can assume that $X$ is
covered by a finite number of convenient neighbourhoods.

For an arbitrary convenient neighbourhood $U_0$, choose an order
of coordinates $u_i$ on it such that $D = \{u_1u_2\cdots u_l =
0\}$.

\begin{prop}\label{clm1}
The zeta-function $\zeta_{\Sig|_{W\setminus Z},\,x}(t)$ at a point
$\,x\in \,U_0\setminus D\,$ is well-defined by the coordinates
$\,u_{l+1}, u_{l+2},\ldots, u_{n+1}\,$ of $x$.
\end{prop}
\begin{prf}
The germ of the manifold $Z$ at the point $x$ is defined by an
equation $\,u_{l+1}^{k_{1,l+1}}\cdots u_{n+1}^{k_{1,n+1}}=0$; in a
neighbourhood of $x$ one has $\,\tilde{F} = \,a\,
u_{l+1}^{k_{2,l+1}}\cdots u_{n+1}^{k_{2,n+1}}$, $\,\Sig =
\,b\,u_{l+1}^{k_{3,l+1}}\cdots u_{n+1}^{k_{3,n+1}}$, where
$\,a(x)\neq 0,\,b(x)\neq 0$, $\,k_{1,j}\in \{0,1\};\, k_{2,j},
k_{3,j} \geq 0.$ The zeta-function $\,\zeta_{\Sig|_{W\setminus
Z},\,x}(t)\,$ is well-defined by the numbers
$\,k_{i,j},\,i=1,2,3,\,j=l+1,\ldots,n+1\,$, which do not depend on
$\,u_1,\ldots,u_l$.
\end{prf}

For a rational function $Q(t)$, we define a set $X_Q = \{x\in
W\cap (\{\Sig =0\}\setminus D) \mid \zeta_{\Sig|_{W\setminus
Z},\,x}(t) = Q(t)\}$. It follows from the proposition above that
for any convenient neighbourhood $U_0$ we have $\,\chi(U_0\cap
X_Q) = 0\,$. Thus for all $Q(t)$ we have $\,\chi(X_Q) = 0\,$  and
$$
\int_{W\cap (\{\Sig =0\}\setminus D)}{\zeta_{\Sig|_{W\setminus
Z},\,x}(t)}^{d\chi} = \prod_{Q}Q^{\chi(X_Q)} = 1.
$$
\end{prf}

\section{Proof of Theorem \ref{thm1}}

Using the Newton diagram $\,\G(F)\,$ of the germ $F$ on
$(\C^{n+1},0)$ one can construct an unimodular simplicial
subdivision $\L$ of the set of covectors with non-negative
coordinates $(\R^{n+1})_+^*$ (see, e.g., \cite{Var}). Consider the
toroidal modification map $\,p:(X_{\L},D) \to (\C^{n+1}, 0)\,$
corresponding to $\L$. Let $\,U\subset \C^{n+1}\,$ be a small
enough ball with the centre at the origin, $\,X = p^{-1}(U)\,$,
$\p=p|_X$. Let $\,Y=\{z_1z_2\cdots z_n = 0\}\subset \C^n_\zz$.
Then $\,S= (Y\times \C_\s)\cup \{\s=0\}\,$ is the union of the
coordinate hyperplanes of $\C^{n+1}$. Since $F$ is non-degenerate
with respect to its Newton diagram $\G(F)$, $\p$ is a resolution
of the germ $S\cup \{F=0\}$ (see, e.g., \cite{Oka}). Finally, $\p$
is an isomorphism outside of $S$, so the resolution $(X,\p)$
satisfies the assumptions of Theorem \ref{thm2}.

Compute the right hand side of (\ref{3}). Let $x\in D\cap W$ be a
point of the $(n-l+1)$-dimensional torus $T_\l$ corresponding to
an $l$-dimensional cone $\l\in \L$. Let $\l$ be generated by
integer covectors $\a_1, \ldots, \a_l$ and let $\l$ lie on the
border of a cone $\l'\in \L$ generated by $\a_1,\ldots,\a_l,\ldots
\a_{n+1}$. Let $(u_1,\ldots, u_{n+1})$ be the coordinate system
corresponding to the set $(\a_1, \ldots, \a_{n+1})$. There exists
a coordinate system $(u_1, \ldots, u_l, w_{l+1}, \ldots, w_{n+1})$
in a neighbourhood $U'$ of the point $x$ such that $\,w_i(x)=0,\,
i = l+1, \ldots, n+1\,$ and $\,\tilde{F}=F\circ
\p=a\,u_1^{k_{1,1}} u_2^{k_{1,2}}\cdots u_l^{k_{1,l}} \cdot
w_{n+1}^{k_{1,n+1}}\,$ (where $\,a(0)\neq 0$). The zero level set
$\,\{\Sig =0\}\,$ is a normal crossing divisor contained in $\{u_1
u_2\cdots u_l = 0\}$. Therefore $\Sig = \s\circ \p = u_1^{k_{2,1}}
u_2^{k_{2,2}}\cdots u_{l}^{k_{2,l}}$. One has: $\,W\cap U' =
\{w_{n+1}=0\}$, $\,(Z\cup \{\Sig =0\})\cap U' = \{u_1 u_2\cdots
u_l = 0\}$, thus $\,\zeta_{\Sig|_{W\setminus Z},x}(t) =
\zeta_{g|_{\{u_i\neq 0,\,i\leq l\}}}(t)$, where $g$ is the germ of
the following function of $n$ variables: $\,g(u_1,\ldots,
u_l,w_{l+1},\ldots, w_n)=u_1^{k_{2,1}} u_2^{k_{2,2}}\cdots
u_l^{k_{2,l}}$.

Assume that one of the exponents $\,k_{2,1}, k_{2,2}\ldots,
k_{2,l}\,$ (say, $k_{2,1}$) is equal to zero. Then $g$ doesn't
depend on $u_1$. We may assume that the monodromy transformation
of its Milnor fibre also doesn't depend on $u_1$. Denote $\,h =
g|_{\{u_1=0\}}\,$. The monodromy transformations of the fibre of
$\,g|_{\{u_2 u_3\cdots u_l \neq 0\}}\,$ and one of $\,h|_{\{u_2
u_3\cdots u_l \neq 0\}}\,$ are homotopy equivalent, so
$\,\zeta_{g|_{\{u_2 u_3\cdots u_l \neq 0\}}}(t)= \zeta_{h|_{\{u_2
u_3\cdots u_l \neq 0\}}}(t).\,$ On the other hand the
multiplicative property of the zeta-function implies that
$\,\zeta_{g|_{\{u_i\neq 0,\,i\leq l\}}}(t)\,\times\,
\zeta_{h|_{\{u_2 u_3\cdots u_l \neq 0\}}}(t)= \zeta_{g|_ {\{u_2
u_3\cdots u_l \neq 0\}}}(t)\,$ and thus $\,\zeta_{g|_{\{u_i\neq
0,\,i\leq l\}}}(t) = 1$.

Now assume that all the exponents $k_{2,1}, k_{2,2}\ldots,
k_{2,l}$ are positive. Then the non-zero fibre of the function $g$
doesn't intersect $\,\{u_1 u_2\ldots u_l = 0\}$, so
$\,\zeta_{g|_{\{u_i\neq 0,\,i\leq l\}}}(t) = \zeta_ {g}(t)$. In
the case $\,l>1\,$ one has $\,\zeta_{g}(t)=1$. In the case
$\,l=1\,$ one has: $\,g = u_1^{k_{2,1}}$, $\,\zeta_g (t) =
1-t^{k_{2,1}}$.

We see that the integrand in (\ref{3}) differs from $1$ only at
points $x$ that lie in strata of dimension $n$. From here on
$\,l=1$. If all the components of $\,\a = \a_1\,$ are positive,
then $\,T_\l \subset D$. Otherwise, $\,T_\l \cap D = \emptyset$.
From here on $\,\a \in P^{\{0,1,\ldots, n\}}\,$ (see the
definitions before Theorem \ref{thm1}).

Using the coordinates $(u_2, \ldots, u_{n+1})$ on the torus
$\,T_{\l} = \{u_1=0\}\,$ we obtain: $T_\l\cap W = \{Q_\a = 0\,\}$,
where for the power series $\,F = \sum F_\kk \s^{k_0}
z_1^{k_1}\cdots z_n^{k_n}\,$ we denote $\,Q_\a=\sum_{\kk\in
\G^{\{0, \ldots, n\}}_\a(F)} F_\kk u_2^{\a_2(\kk)}
u_3^{\a_3(\kk)}\cdots u_{n+1}^{\a_{n+1}(\kk)}.\,$ So $\,\,T_\l\cap
W\,$ is the zero level set of the Laurent polynomial $Q_\a$. Using
results of \cite{Khov1}, \cite{Khov2} we obtain: $\chi(T_\l\cap W)
= (-1)^{n-1}n!\,V_n(\D(Q_\a))\,$, where $\D(\cdot)$ denotes the
Newton polyhedron. Since the polyhedra $\D(Q_\a)$ and $\,\G_\a =
\G^{\{0,1,\ldots,n\}}_\a (F)\,$ are isomorphic as subsets of
integer lattices, their volumes are equal: $\,V_n(\D(Q_\a)) =
V_n(\G_\a)$. In a neighbourhood of a point $\,x\in T_\l\cap W\,$
one has $\Sig = a\,u_1^{\a(\partial/\partial k_0)}$, where
$a(x)\neq 0$. Therefore $\zeta_{\Sig|_{W\setminus Z},\,x}(t) =
1-t^{\a(\partial/\partial k_0)}$. Thus one has:
\begin{equation}
\label{5} \int_{T_\l\cap W} {\zeta_{\Sig|_{W\setminus
Z},\,x}(t)}^{d\chi} = (1-t^{\a(\frac{\partial}{\partial
k_0})})^{\chi(T_\l\cap W)} = (1-t^{\a(\frac{\partial}{\partial
k_0})})^{(-1)^{n-1}n!\,V_n(\G_\a)}.
\end{equation}

Multiplying (\ref{5}) for all strata $T_\l \subset D$ of dimension
$n$ we get (\ref{1}).

\begin{acknl}
I am grateful to my advisor S. M. Gusein-Zade, who motivated me to
make this study, for constant attention and support.
\end{acknl}

\noindent E-mail: gusev@mccme.ru

\end{document}